\input amstex
\input amsppt.sty   
\hsize 26pc
\vsize 42pc
\magnification=\magstep1
\def\nmb#1#2{#2}         
\def\totoc{}             
\def\idx{}               
\def\ign#1{}             

\define\X{\frak X}
\define\g{\frak g}
\redefine\o{\circ}
\define\al{\alpha}
\define\be{\beta}
\define\ga{\gamma}

\define\ze{\zeta}
\define\et{\eta}

\define\ka{\kappa}
\define\la{\lambda}
\define\rh{\rho}
\define\si{\sigma}
\define\ta{\tau}
\define\ph{\varphi}

\define\ps{\psi}
\define\om{\omega}
\define\Ga{\Gamma}

\define\Ph{\Phi}

\define\Om{\Omega}
\redefine\i{^{-1}}
\define\row#1#2#3{#1_{#2},\ldots,#1_{#3}}
\define\x{\times}
\define\uet{\underline{\et}}
\define\Fl{\operatorname{Fl}}
\define\Pt{\operatorname{Pt}}
\def\today{\ifcase\month\or
 January\or February\or March\or April\or May\or June\or
 July\or August\or September\or October\or November\or December\fi
 \space\number\day, \number\year}
\topmatter
\title The relation between systems \\ 
and associated bundles  \endtitle
\author  Peter W. Michor  \endauthor
\affil
Institut f\"ur Mathematik, Universit\"at Wien,\\
Strudlhofgasse 4, A-1090 Wien, Austria.
\endaffil
\address{P. Michor: Institut f\"ur Mathematik, Universit\"at Wien,
Strudlhofgasse 4, A-1090 Wien, Austria.}\endaddress
\email MICHOR\@AWIRAP.bitnet \endemail
\date{\today}\enddate
\thanks{This work was done during a visit of the author at
Dipartimento di Matematica Applicata, Universit\`a di Firenze
Via San Marta, 3, I-50139 Firenze, Italy, during July 1990, which was
supported by a grant from Gruppo 
Nazionale Fisica Matematica, CNR, Italy.}\endthanks
\abstract{It is shown that a strong system of vector fields on a fiber 
bundle in the sense of \cite{Mo} is induced from a principal fiber 
bundle if and only if each vertical vector field of the system is 
complete.}\endabstract
\endtopmatter
\leftheadtext{\smc Peter W. Michor}
\rightheadtext{\smc Relation between systems and associated bundles}
\document

\heading Table of contents \endheading
\noindent 1. Introduction \leaders \hbox to 1em{\hss .\hss }\hfill 
  {\eightrm 1}\par 
\noindent 2. Systems \leaders \hbox to 1em{\hss .\hss }\hfill {\eightrm 2}\par 
\noindent 3. Properties of complete strong systems \leaders 
  \hbox to 1em{\hss .\hss }\hfill {\eightrm 6}\par 
\noindent 4. The universal connection \leaders \hbox to 1em{\hss .\hss }
  \hfill {\eightrm 12}\par 

\heading\totoc\nmb0{1}. Introduction \endheading

The notion of systems of vector fields and systems of connections for 
fibered manifolds
were introduced by Marco Modugno as a generalisation of principal 
connections and as a means to give a lucid and easy construction of the 
universal connection on the bundle of connections, which for 
principal bundles is due to \cite{Garcia, 1977}. This is a special 
case of the usual notion of a system as treated for example in 
\cite{Gauthier, 1984}

In this paper we prove as the main result (theorem \nmb!{3.6})
that any such system for a fiber bundle which 
is strong in the sense of Modugno and has the further property that 
all vertical vector fields of the system are complete, is in fact 
an induced system on an associated bundle for a principal bundle. 
The structure group of the principal bundle is the holonomy group of 
the universal connection of the system. In \nmb!{2.7} we show that 
the converse is true and we describe some simple examples of 
non-complete systems also.

We use heavily the concepts and techniques of \cite{Mi}. These can be 
found also with more details and more complete proofs in 
\cite{Michor, 1991}.

I want to thank Marco Modugno for his hospitality, for asking the 
question answered in this paper, and for lots of discussions.

\heading\totoc\nmb0{2}. Systems \endheading

\subheading{\nmb.{2.1}} Let $(E,p,M,S)$ be a smooth finite 
dimensional fibre bundle with base $M$ and standard fibre $S$.
By a \idx{\it system of vector fields} on $E$ we mean a pair 
$(H,\et)$, where
\roster
\item $H=(H,q_H,M)$ is a vector bundle over $M$.
\item $\et$ is a mapping factoring as in the following diagram
$$\CD
H\x_ME   @>\et>>  TE\\
@Vpr_1VV          @VVTpV\\
H        @>>\uet> TM.
\endCD$$
Furthermore it is supposed to be fibre respecting and fiber linear over 
$E$. Then it turns out that $\underline{\et}:H\to TM$ is fibre 
respecting and fiber linear over $M$. Also we suppose that 
$\underline{\et}$ is fiberwise surjective.
\endroster
In \cite{Mo} this is called a linear, horizontally complete, and 
projectable system of vector fields.

\subheading{\nmb.{2.2}} A system of vector fields is called 
\idx{\it monic} if the associated mapping 
$\check \et: H\to \bigcup_{x\in M}C^\infty(TE|E_x)$
is injective.

It is called \idx{\it involutive}, if the push forward of the associated 
mapping acting on sections
$\check \et_*: C^\infty(H)\to C^\infty(TE)=\X(E)$ has as image a Lie 
sub algebra of the algebra of vector fields on $E$. Note that this is 
not involutivity of some sub bundle of $TE$, since all vector fields 
in the image of $\check \et_*$ are "rigid" along the fibres of $E$.

A system is called \idx{\it canonical} if there exist an open cover 
$(U_\al)$ of $M$, a fiber bundle atlas
$(U_\al,\ps_\al:E|U_\al\to U_\al\x S)$ of $(E,p,M)$, and a vector 
bundle atlas $(U_\al, \ph_\al:H|U_\al\to TU_\al\x V)$  of 
$(H,q_H,M)$, such that
$$T\ps_\al.\et.(\ph_\al\i(\xi_x,v), 
     \ps_\al\i(x,s))=(\xi_x,\et^\al(v)(s)),$$
where $\et^\al: V\to \X(S)$ is a linear mapping into the space of 
vector fields on the standard fiber $S$. So it is required that the 
the mapping $\et^\al$ does not depend on the foot point $x\in U_\al$.
These data will be called \idx{\it canonical atlases} for the system.

A system that is monic, involutive, and canonical is called a 
\idx{\it strong system}, see \cite{Mo}.

\subheading{\nmb.{2.3}} Let $(H,\et)$ be a system of vector fields on 
the bundle $E$. Then the kernel of the vector bundle homomorphism 
$\uet: H\to TM$ is a sub vector bundle $A$ of $H$. Thus we have the 
following diagram
$$\CD
  @.     @.    H\x_ME   @>\et>>  TE    @.\\
@.     @.      @Vpr_1VV          @VVTpV   @.\\
0 @>>> A @>>i> H        @>>\uet> TM    @>>> 0,
\endCD$$
where the lower line is an exact sequence of vector bundles.
For $a_x\in A_x$ the vector field $\check\et(a_x)\in \X(E_x)$.

We say that the system $(H,\et)$ is \idx{\it complete} if and only if 
each vector field $\check\et(a_x)\in \X(E_x)$ is a complete vector 
field on the fiber $E_x$; so its flow should exist for all time.

\subheading{\nmb.{2.4}. Remark} The exact sequence 
$0 \to A \to H \to TM$ of a monic and involutive system is also 
called a \idx{\it Lie algebroid}, see e\. g\. \cite{Mackenzie, 1987, 
p\. 100}, or an \idx{\it abstract Atiyah sequence}, see 
\idx{\it Almeida-Molino, 1985}; one forgets the bundle $E$ on which 
the sections of $H$ induce projectable vector fields. In this paper 
we will concentrate on $E$.

\subheading{\nmb.{2.5}} 
Let $(H,\et)$ be an involutive monic system of vector fields on 
the bundle $(E,p,M,S)$. We consider the exact sequence 
$0\to A\to H \to TM \to 0$ of vector bundles from \nmb!{2.3} and the 
induced exact sequence of push forwards on the respective spaces of 
sections
$$0 @>>> C^\infty(A) @>i_*>> C^\infty(H) @>\uet_*>> \X(M) \to 0.$$
We have also the induced push forward mapping
$$\check\et_*:C^\infty(H)\to \X(E)$$
which is injective since the system is monic.
The image of $\check\et_*$ is closed under the Lie bracket, so there 
is an induced bracket 
$$[\quad,\quad]^H:C^\infty(H)\x C^\infty(H) \to C^\infty(H)$$ 
which is a bilinear differential operator of total degree 1.

Now for sections $a_1$, $a_2\in C^\infty(A)$ and a function 
$f\in C^\infty(M,\Bbb R)$ we have 
$[\check\et_*(f.a_1), \check\et_*(a_2)] = f.[\check\et_*(f.a_1), 
\check\et_*(a_2)]$, which is again vertical, 
since the vector fields $\check\et_*(a_i)\in \X(E)$ are vertical and 
$f$ is constant along the fibres.
Thus the induced bracket $[\quad,\quad]: 
C^\infty(A)\x C^\infty(A)\to C^\infty(A)$ is of order 0 und is thus a 
push forward by a smooth fiberwise Lie bracket 
$[\quad,\quad]^A:A\x_MA\to A$. 
Note that the isomorphism type of the Lie algebra 
$(A_x,[\quad,\quad]^A_x)$ need not be locally constant, if the Lie 
algebra is not rigid, for example.

Let us now assume furthermore that the monic involutive system is 
also canonical (see \nmb!{2.2}) and let $(U_\al,\ps_\al)$ and $(U_\al, 
\ph_\al)$ be canonical atlases for this system as spelled out in 
\nmb!{2.2}. We want to express the bracket $[h_1,h_2]^H$ for $h_1$, 
$h_2\in C^\infty(H)$ in terms of the canonical atlases.
We have $\ph_\al(h_i)(x)=(X_{h_1}(x), v_{h_i}(x))\in TU_\al\x V$, and 
$h_i$ has values in the sub bundle $A$ if and only if the vector 
field $X_{h_i}$ is zero.
We have then
$$\align
(T&\ps_\al\o \check\et_*(h))(\ps_\al\i(x,s))
     = (T\ps_\al\o\et)(\ph_\al\i(X_h(x),v_h(x)),\ps_\al\i(x,s))\\
&=(X_h(x),\et^\al(v_h(x))(s))\\
(T&\ps_\al\o [\check\et_*(h_1),\check\et_*(h_2)]^{\X(E)})(\ps_\al\i(x,s))=\\
&=[(X_{h_1},\et^\al\o v_{h_1}), 
     (X_{h_2},\et^\al\o v_{h_2})]^{\X(U_\al\x S)}(x,s)\\
&=\Bigl([X_{h_1},X_{h_2}]^{\X(U_\al)}(x), 
    [\et^\al(v_{h_1}(x)),\et^\al(v_{h_2}(x))]^{\X(S)}(s)\\
&\qquad     +d(\et^\al\o v(h_2))(x)(X_{h_1}(x))(s)
     - d(\et^\al\o v(h_1))(x)(X_{h_2}(x))(s)\Bigr).
\endalign$$
If $a_1$, $a_2\in C^\infty(A)$ then we get
$$\align
(T\ps_\al&\o [\check\et_*(a_1),\check\et_*(a_2)]^{\X(E)})(\ps_\al\i(x,s))=\\
&= \Bigl(0,[\et^\al(v_{a_1}(x)),\et^\al(v_{a_2}(x))]^{\X(S)}(s)\Bigr)\\
&=: \et^\al([v_{a_1}(x),v_{a_2}(x)]^V)(s).
\endalign$$
So the canonical atlases for a canonical system restrict to a vector 
bundle atlas for the Lie algebra bundle $(A,[\quad,\quad]^A)$ in 
which the Lie algebra structure is locally trivial, thus constant along 
connected components of $M$. To simplify notation we assume that it 
is constant, isomorphic to $V,[\quad,\quad]^V$. 

\subheading{\nmb.{2.6}. Connections for a system}
Let $(H,\et)$ be a system of vector fields on 
the bundle $E$. We consider a vector bundle homomorphism 
$\si:TM\to H$ which splits the exact sequence 
$0\to A\to H \to TM \to 0.$
Then $\si$ defines a horizontal lifting $C_\si:TM\x_ME\to TE$ by the 
prescription 
$$C_\si(\xi_x, u_x):= \et(\si(\xi_x),u_x)\in TE.$$
So $C_\si$ is linear over $E$ and is a right inverse to $(Tp,\pi_E): 
TE \to TM\x_ME$. By \cite{Mi, 1.1} $C_\si$ specifies a connection for 
the fiber bundle $E$. We call all connections obtained in this way 
\idx{\it connections respecting the system} $H$ or just 
\idx{\it $H$-connections}.

Let us suppose now for the moment that the system $H$ is canonical 
and let $(U_\al,\ps_\al)$, $(U_\al, \ph_\al)$ be canonical atlases for 
the system $H$ as required in \nmb!{2.2}.  The splitting 
$\si:TM\to H$ can then written as 
$\ph_\al(\si(\xi_x))=(\xi_x,\si^\al(\xi_x))$, where 
$\si^\al\in\Om^1(U_\al;V)$ is a one form on $U_\al$ with values in 
the vertical part $V$ of the standard fiber of $H$.

The space of all vector bundle splittings of the exact sequence 
$0\to A\to H \to TM \to 0$ parametrizes thus  the space of all 
connections of the fiber bundle $E$ which respect the system 
$(H,\et)$. These splittings are exactly the sections of the affine 
bundle 
$$C(H):=\{s_x\in L(T_xM,H_x): \uet_x\o s_x=Id_{T_xM},x\in M\}.$$ 
The modelling bundle of that affine bundle is $L(TM,A)=T^*M\otimes A$.

\subheading{\nmb.{2.7}. Associated systems}
Let $(P,M,p,G)$ be a principal fiber bundle with structure group $G$, 
and let $\ell: G\x S\to S$ be a smooth left action on a smooth 
manifold. Then we have the \idx{\it associated fiber bundle}
$P[S]=P[S,\ell]=(P\x S)/G$. On the principal bundle $P$ there is the 
strong system of all projectable $G$-equivariant vector fields 
$(TP/G,\et_P)$, whose exact sequence in the sense of \nmb!{2.3} is 
given by
$$0 \to VP/G=P[\g, Ad] \to TP/G \to TM \to 0.$$
The sections of $TP/G$ correspond to the infinitesimal automorphisms of the 
principal bundle. The vertical sections correspond to the 
infinitesimal gauge transformations. 

The strong system $(TP/G,\et_P)$ thus 
induces a system $(TP/G,\et_{P[S]})$ on the 
associated bundle $P[S]$ which is monic if and only if the action 
$\ell$ is infinitesimally effective, i\. e\. the fundamental vector 
field mapping $\ze:\g\to \X(S)$ is injective. By looking at a 
principal bundle atlas and the induced associated atlas (see 
\cite{Mi, section 2}) one easily sees that these systems are 
canonical and complete. Also it is easily checked, that an arbitrary 
system $(H,\et)$ on the associated bundle $P[S]$ is isomorphic to the 
induced system if and only if each $H$-connection is induced from a 
principal connection; by using \cite{Mi, 2.5} one may recognize 
these induced connections.

If we take a suitable open subbundle $E$ of the associated bundle 
$P[S]$ we obtain by restriction a (strong) system $(TP/G,\et_E)$ which in 
general is not complete. 

\heading\totoc\nmb0{3}. Properties of complete strong systems \endheading

\proclaim{\nmb.{3.1}. Theorem} 
Let $(H,\et)$ be a complete strong system of vector fields on 
the bundle $(E,p,M,S)$.

Then each connection $C_\si$ respecting the system $H$ for any 
splitting of the exact sequence $0\to A\to H \to TM \to 0$ is 
complete in the sense of \cite{Mi, 1.6}: its parallel transport 
exists globally.
\endproclaim

\demo{Proof}
Let $c:[0,1] \to  M$ be a smooth curve. 
We have to show that for each $u_0\in E_{c(0)}$ there exists a smooth 
curve $\operatorname{Pt}(c,t,u_0)$ in $E$ which covers $c(t)$, 
is horizontal, has initial value $u_0$, and is defined for all 
$t\in [0,1]$. We refer to \cite{Mi, theorem 1.5} for the local 
existence and general properties of parallel transport.

Let $(U_\al,\ps_\al)$, $(U_\al, \ph_\al)$ be canonical atlases for 
the system $H$ as required in \nmb!{2.2}.
We choose a partition 
$0 = t_0 < t_1 < \cdots < t_k = 1$ such that $c([t_i,t_{i+1}])
\subset U_{\al_i}$ for suitable $\al_i$. It suffices to show
that $\operatorname{Pt}(c(t_i + \quad),t,u_{c(t_i)})$ exists for
all $0\leq t \leq t_{i+1} - t_i$ and all $u_{c(t_i)}\in E_{c(t_i)}$, 
for all $i$ --- then we may piece
them together. So we may assume that $c:[0,1] \to  U_\al$ for some
$\al$.

By \cite{Mi, third proof of 1.5} we have in $U_\al\x S$
$$\ps_\al(\operatorname{Pt}(c,t,\ps_\al\i(c(0),s))) = 
     (c(t),\ga(s,t)),$$
where $\ga(s,t)$ is the evolution line (integral curve) of the time 
dependent vector field $\Ga^\al(\frac d{dt}c(t))$ on $S$, where 
$\Ga^\al\in \Om^1(U_\al,\X(S))$ is the Christoffel form for the 
connection $C_\si$ in the fiber bundle chart $(U_\al,\ps_\al)$, see 
\cite{Mi, 1.4}, from where we use now the defining equation for 
$\Ga^\al$ to compute as follows, where $\Ph_\si$ is the projection 
onto the vertical bundle $VE$ along the horizontal bundle 
$C_\si(TM\x_ME)$:
$$\align
(0_x,&\Ga^\al(\xi_x,s)) = - T(\ps_\al).\Ph_\si.T(\ps_\al)\i (\xi_x,0_s)\\
&= - T(\ps_\al).(T(\ps_\al)\i (\xi_x,0_s)-C_\si(\xi_x,\ps_\al\i(x,s)))\\
&= -(\xi_x,0_s) + T(\ps_\al).\et.(\si(\xi_x),\ps_\al\i(x,s))\\
&= -(\xi_x,0_s) 
     + T(\ps_\al).\et.(\ph_\al\i(\xi_x,\si^\al(\xi_x)),\ps_\al\i(x,s))\\
&= -(\xi_x,0_s) + (\xi_x,\et^\al(\si^\al(\xi_x))(s))\\
&= (0_x,\et^\al(\si^\al(\xi_x))(s)).
\endalign$$ 

Since the system $H$ is complete by assumption we have 
$$T(\ps_\al).\et.(\ph_\al\i(0_x,v),\ps_\al\i(x,s))= (0_x,\et^\al(v)(s))$$
and $\et^\al(v)\in\X(S)$ is a complete vector field for each 
$v\in V$. So $\et^\al:(V,[\quad,\quad]^V)\to \X(S)$ is a 
homomorphism of Lie algebras whose image consists of complete vector 
fields. By the theorem of \cite{Palais, 1957} there is a simply 
connected Lie group $G_\al$ with Lie algebra $V$ 
and a right action $r_\al:S\x G_\al\to S$ 
of $G_\al$ on $S$ such that $\et^\al$ is the fundamental vector field 
mapping for this action: $\et^\al(v)(s)=T_e(r_\al(s,\quad))v$.

 From the computation above we have
$\Ga^\al(\frac d{dt}c(t))=\et^\al(\si^\al(\frac d{dt}c(t)))$.
Let us choose a left invariant Riemannian metric on the Lie group 
$G_\al$. It is then a complete Riemannian metric, and the left 
invariant vector fields $L(v)$ generated by the $v\in V$ are all 
bounded with respect to this metric. Since $[0,1]$ is compact,
$L(\si^\al(\frac d{dt}c(t)))$ is a time dependent vector field which 
is bounded for the complete metric. Thus there exists the global 
evolution curve $t\mapsto g_\al(t)\in G_\al$ for $t\in[0,1]$, 
uniquely given by 
$$\cases \tfrac d{dt} g_\al(t) = 
     T\la_{g_\al(t)}.\si^\al(\frac d{dt}c(t))\\
g_\al(0)=e,
\endcases$$
where $\la_g$ is left translation by $g\in G_\al$.
But then we have
$$\align
\frac d{dt} r_\al(s, g_\al(t)) &= T(r_\al(s,\quad)) 
     \frac d{dt}g_\al(t) \\
&= T(r_\al(s,\quad)). T_e(\la_{g_\al(t)}).\si^\al(\frac d{dt}c(t))\\
&= T_e(r_\al(r_\al(s,g_\al(t)),\quad)). \si^\al(\frac d{dt}c(t))\\
&= \et^\al(\si^\al(\frac d{dt}c(t)))(r_\al(s,g_\al(t))),\\
r_\al(s,g_\al(0))&= r_\al(s,e)=s.
\endalign$$
Thus $r_\al(s,g_\al(t))=\ga(s,t)$, the looked for global evolution 
curve for for the time dependent vector field
$\Ga^\al(\frac d{dt}c(t))=\et^\al(\si^\al(\frac d{dt}c(t)))$.
\qed\enddemo

\subheading{\nmb.{3.2}. Curvature}
Let $(H,\et)$ be a complete strong system of vector fields on 
the bundle $(E,p,M,S)$.
We want to compute the curvature $R$ of a $H$-respecting connection
$C=C_\si$ for a splitting $\si:TM\to H$ in canonical coordinates.

 From \cite{Mi, 1.4} we have 
$$((\ps_\al\i)^*R)((X_1,Y_1),(X_2,Y_2)) = d\Ga^\al(X_1,X_2) + 
[\Ga^\al(X_1),\Ga^\al(X_2)]^{\X(S)}.$$
 From the proof of \nmb!{3.1} we have 
$\Ga^\al=\et^\al\o \si^\al\in \Om^1(U_\al,\X(S))$, 
thus we may compute 
$$\align
((\ps_\al\i)^*R)&((X_1,Y_1),(X_2,Y_2)) = d\Ga^\al(X_1,X_2) 
     + [\Ga^\al(X_1),\Ga^\al(X_2)]^{\X(S)}\\
&= d(\et^\al\o\si^\al)(X_1,X_2) 
     + [(\et^\al\o\si^\al)(X_1),(\et^\al\o\si^\al)(X_2)]^{\X(S)}\\
&= \et^\al\Bigl(d\si^\al(X_1,X_2) + [\si^\al(X_1),\si^\al(X_2)]^V\Bigr).
\endalign$$

\subheading{\nmb.{3.3}. The holonomy Lie algebra} 
Let $(H,\et)$ be a complete strong system. The holonomy Lie 
algebra of any 
(complete by \nmb!{3.1}) $H$-connection $C_\si$ is given as follows (see 
\cite{Mi}, 3.2): 

Let $M$ be connected. Choose $x_0\in M$, a base point, and identify the 
standard fiber $S$ with $E_{x_0}$. For $x\in M$ and $X_x$, 
$Y_x\in T_xM$ we consider the horizontal lifts $C_\si(X_x)$ and 
$C_\si(Y_x)$ which are vector fields on $E$ along $E_x$. Then the 
curvature applied to these fields is vertical,  
$R(C_\si(X_x),C_\si(Y_x))\in\X(E_x)$. Now we choose a piecewise 
smooth curve $c$ in $M$ from $x_0$ to $x$ and consider the pullback 
under the parallel transport 
$$\operatorname{Pt}(c,1,\quad)^*R(C_\si(X_x),C_\si(Y_x))
     \in\X(E_{x_0})=\X(S).$$ 
The closed linear span of all these vector fields in $\X(S)$ with 
respect to the compact $C^\infty$-topology is called the \idx{\it holonomy 
Lie algebra} $\operatorname{hol}(C_\si,x_0)$ of the connection 
$C_\si$, centered at $x_0$.

\proclaim{\nmb.{3.4}. Lemma} The holonomy Lie algebra  
$\operatorname{hol}(C_\si,x_0)$ is a sub Lie algebra of 
$\et_{x_0}(A_{x_0})\subset \X(E_{x_0})$ and is thus finite 
dimensional.
\endproclaim

\demo{Proof}
Using \nmb!{3.2} and the proof of \nmb!{3.1} we get in turn
$$\align
T(\ps_\al)&.R(C_\si(X_x,\ps_\al\i(x,s)),C_\si(Y_x,\ps_\al\i(x,s))) \\
&= ((\ps_\al\i)^*R)(T(\ps_\al).\et.
     (\ph_\al\i(X_x,\si^\al(X_x)),\ps_\al\i(x,s))),\ldots) \\
&= ((\ps_\al\i)^*R)((X_x,(\et^\al\o\si^\al)(X_x)(s)),
     (Y_x,(\et^\al\o\si^\al)(Y_x)(s))) \\
&= \et^\al\Bigl(d\si^\al(X_x,Y_x) + 
     [\si^\al(X_x),\si^\al(Y_x)]^V\Bigr)(s).
\endalign$$
Thus $R(C_\si(X_x),C_\si(Y_x))\in \et_x(A_x)$. Next we prove that 
pull back via parallel transport does not move out of $\et(A)$.
 From the proof of theorem \nmb!{3.1} we know that for a smooth curve
$c$ in $U_\al$ we have 
$$\ps_\al(\operatorname{Pt}(c,t,\ps_\al\i(c(0),s))) = 
     (c(t),r_\al(s,g_\al(t))),$$
where $g_\al(t)$ is a globally defined curve in the Lie group $G_\al$ 
and where $r_\al:S\x G_\al\to S$ is a right action such that 
$\et_\al:V\to\X(S)$ is the fundamental vector field mapping.
But then we have 
$$\align
&(\ps_\al\i)^* \operatorname{Pt}(c,t)^*\check\et(\ph_\al\i(0_{c(t)},v)) = \\
&= T\ps_\al\o T\operatorname{Pt}(c,t)\i\o 
     \check\et(\ph_\al\i(0_{c(t)},v)) \o \operatorname{Pt}(c,t) 
     \o \ps_\al\i(c(0),\quad) \\
&= T(\ps_\al\o \operatorname{Pt}(c,t)\o \ps_\al\i)
     (0_{c(t)},\et^\al(v)\o \operatorname{pr}_2\o \ps_\al\o 
     \operatorname{Pt}(c,t)\o\ps_\al\i(c(0),\quad))  \\
&= (0_{c(0)}, T(r_\al^{g_\al(t)\i})\o \et^\al(v)\o r_\al^{g_\al(t)}) =\\
&= (0_{c(0)}, (r_\al^{g_\al(t)})^*\et^\al(v)) = \\
&= (0_{c(0)}, \et^\al(\operatorname{Ad}(g_\al(t))v)), 
\endalign$$
by well known properties of right Lie group actions.
This implies the desired result.
\qed\enddemo

\subheading{\nmb.{3.5}. Holonomy groups} Let
(E,p,M,S) be a fibre 
bundle with a complete connection $\Ph$, and let us assume that
$M$ is connected. We choose a fixed base point $x_0\in M$
and we identify $E_{x_0}$ with the standard fiber $S$. For each
closed piecewise smooth curve $c:[0,1]\to M$ through $x_0$ the
parallel transport $\Pt(c,\quad,1) =: \Pt(c,1)$ (pieced together
over the smooth parts of $c$) is a diffeomorphism of $S$. All
these diffeomorphisms form together the group
$\operatorname{Hol}(\Ph,x_0)$, the \idx{\it holonomy group} of $\Ph$
at $x_0$, a subgroup of the diffeomorphism group
$\operatorname{Diff}(S)$. If we consider only those piecewise
smooth curves which are homotopic to zero, we get a subgroup
$\operatorname{Hol}_0(\Ph,x_0)$, called the \idx{\it restricted holonomy group} of
the connection $\Ph$ at $x_0$.

\proclaim{\nmb.{3.6}. Theorem} 
Let $(H,\et)$ be a complete strong system of vector fields on 
the bundle $(E,p,M,S)$. Let $M$ be connected.

Then there is a principal bundle $(P,p,M,G)$ with finite
dimensional structure group $G$ and a smooth action of $G$ 
on $S$ such that the following statements hold.
\roster
\item The Lie algebra $\g$ of $G$ is isomorphic to 
     the standard Lie algebra $(V,[\quad,\quad])$ of the 
     Lie algebra bundle $(A,p,M)$.
\item The fibre bundle $E$ is
     isomorphic to the associated bundle $P[S]$.
\item The system $(H,p,M)$ is induced from the system $TP/G$ of 
     $G$-invariant projectable vector fields on $P$.
\item Any connection respecting $H$ is induced from a principal 
     connection on $P$.
\endroster
\endproclaim

\demo{Proof} 
Let us again identify $E_{x_0}$ and $S$, and also $A_{x_0}$ and the 
standard Lie algebra of the Lie algebra bundle $(A,p,M)$. 
Then $\et_{x_0}:A_{x_0}\to \X(E_{x_0})$ is a Lie algebra homomorphism 
whose image consists of complete vector fields, since the system is 
complete. There exists a Lie group $G_0$ with Lie algebra $A_{x_0}$ 
and a effective smooth left action 
$\ell:G_0\x S\to S$ such that $-\et_{x_0}$ is 
the fundamental vector field mapping  for it (which is a Lie algebra anti 
homomorphism for left actions). We call $\g$ the image of 
$\et_{x_0}$. This is then a finite dimensional sub Lie algebra of 
$\X(S)$ which is anti isomorphic to the Lie algebra of $G_0$. We view 
$G_0$ as a finite dimensional subgroup of the group of all 
diffeomorphisms of $S=E_{x_0}$.
For the rest of the proof we choose an $H$-connection $C$, 
given by some splitting $\si:TM\to H$, which we fix from 
now on. 

\remark{Claim 1} $G_0$ contains $\operatorname{Hol}_0(\Ph,x_0)$, the
restricted holonomy group. \newline
Let $f \in \operatorname{Hol}_0(\Ph,x_0)$, then $f=\Pt(c,1)$ for
a piecewise smooth closed curve $c$ through $x_0$, which is
nullhomotopic. Since the parallel transport is essentially
invariant under reparametrisation, see \cite{Mi, 1.5.3}, we can replace $c$
by $c\o g$, where $g$ is smooth and flat at each corner of $c$.
So we may assume that $c$ itself is smooth. Since $c$ is
homotopic to zero, by approximation we may assume that there is
a smooth homotopy $H: \Bbb R^2\to M$ with $H_1|[0,1] = c$ and 
$H_0|[0,1] = x_0$. Then $f_t:= \Pt(H_t,1)$ is a curve in
$\operatorname{Hol}_0(\Ph,x_0)$ which is smooth as a mapping
$\Bbb R\x S\to S$.
\endremark

\remark{Claim 2} $(\tfrac d{dt} f_t)\o f_t\i =: Z_t$ is in 
$\g$ for all $t$. 
\endremark
To prove claim 2 we consider the pullback bundle $H^*E\to \Bbb
R^2$ with the induced connection $H^\Ph$. It is sufficient to
prove claim 2 there. Let $X=\tfrac d{ds}$ and $Y=\tfrac d{dt}$
be the constant vector fields on $\Bbb R^2$, so $[X,Y]=0$. Then 
$\Pt(c,s) = \Fl^{CX}_s|(H^*E)_{(1,0)}$ and so on. We put 
$$f_{t,s}= \Fl^{CX}_{-s}\o\Fl^{CY}_{-t}\o\Fl^{CX}_s\o\Fl^{CY}_t:S\to S,$$
so $f_{t,1}=f_t$. Then we have in the vector space $\X(S)$
$$\multline 
(\tfrac d{dt} f_{t,s})\o f_{t,s}\i = -(\Fl^{CX}_s)^*CY 
    + (\Fl^{CX}_s)^*(\Fl^{CY}_t)^*(\Fl^{CX}_{-s})^*CY, \\ \allowdisplaybreak
\shoveleft{(\tfrac d{dt} f_{t,s})\o f_{t,s}\i 
    = \int_0^s\tfrac d{ds}\left((\tfrac d{dt} f_{t,s})\o f_{t,s}\i\right)ds}\\
\quad = \int_0^s\left(-(\Fl^{CX}_s)^*[CX,CY] 
    + (\Fl^{CX}_s)^*[CX,(\Fl^{CY}_t)^*(\Fl^{CX}_{-s})^*CY]\right. \\
    \left.- (\Fl^{CX}_s)^*(\Fl^{CY}_t)^*(\Fl^{CX}_{-s})^*[CX,CY]\right)\;ds.
\endmultline$$
Since $[X,Y]=0$ we have $[CX,CY] = \Ph[CX,CY] = R(CX,CY)$ and
$$\multline (\Fl^{CX}_t)^*CY = C\left((\Fl^{X}_t)^*Y\right) 
     + \Ph\left((\Fl^{CX}_t)^*CY\right) \\
= CY + \int_0^t\tfrac d{dt}\Ph(\Fl^{CX}_t)^*CY\;dt 
= CY + \int_0^t\Ph(\Fl^{CX}_t)^*[CX,CY]\;dt \\
= CY + \int_0^t\Ph(\Fl^{CX}_t)^*R(CX,CY)\;dt. 
\endmultline$$
Thus all parts of the integrand above are in $\g$ and so 
$(\tfrac d{dt} f_{t,s})\o f_{t,s}\i$ is in $\g$ for all $t$
and claim 2 follows.

Now claim 1 can be shown as follows. There is a unique smooth
curve $g(t)$ in $G_0$ satisfying $T_e(\rho_{g(t)})Z_t = Z_t.g(t)
= \tfrac d{dt} g(t)$
and $g(0) = e$ where $\rh_g$ denotes right translation by $g$ in $G$. 
Via the action of $G_0$ on $S$ the curve $g(t)$
is a curve of diffeomorphisms on $S$, generated by the time
dependent vector field $Z_t$, so $g(t) = f_t$ and $f =f_1$ is in
$G_0$. So we get $\operatorname{Hol}_0(\Ph,x_0)\subseteq G_0$.

\remark{Step 3} Now let $G$ be the subgroup of the group of all 
diffeomorphisms of $S$ which is generated by the full holonomy group 
$\operatorname{Hol}(\Ph,x_0)$ and by $G_0$. We make $G$ into a
Lie group by taking $G_0$
as its connected component of the identity. This is possible:
$G/G_0$ is a countable group, since the
fundamental group $\pi_1(M)$ is countable (by Morse Theory $M$
is homotopy equivalent to a countable CW-complex).
\endremark

\remark{Step 4} Construction of a cocycle of transition functions
with values in $G$. Let $(U_\al,u_\al:U_\al\to \Bbb R^m)$ be a
locally finite smooth atlas for $M$ such that each
$u_\al:U_\al\to \Bbb R^m)$ is surjective. Put $x_\al :=
u_\al\i(0)$ and choose smooth curves $c_\al:[0,1]\to M$ with
$c_\al(0)=x_0$ and $c_\al(1)=x_\al$. For each $x\in U_\al$ let
$c_\al^x:[0,1]\to M$ be the smooth curve $t\mapsto u_\al\i(t.u_\al(x))$,
then $c_\al^x$ connects $x_\al$ and $x$ and the mapping
$(x,t)\mapsto c_\al^x(t)$ is smooth $U_\al\x[0,1]\to M$.
Now we define a fibre bundle atlas 
$(U_\al,\ps_\al:E|U_\al\to U_\al\x S)$ by 
$\ps_\al\i(x,s)= \Pt(c_\al^x,1)\,\Pt(c_\al,1)\,s$. Then $\ps_\al$
is smooth since $\Pt(c_\al^x,1) = \Fl^{CX_x}_1$ for a local
vector field $X_x$ depending smoothly on $x$. Let us investigate
the transition functions. 
$$\align \ps_\be\ps_\al\i(x,s) &= \left(x,
    \Pt(c_\al,1)\i \Pt(c_\al^x,1)\i \Pt(c_\be^x,1)\Pt(c_\be,1)\,s\right)\\
&= \left(x,\Pt(c_\be.c_\be^x.(c_\al^x)\i.(c_\al)\i,4)\,s\right)\\
&=: (x,\ps_{\be\al}(x)\,s),\text{ where }\ps_{\be\al}:U_{\be\al}\to G.
\endalign$$
Clearly $\ps_{\be\al}: U_{\be\al}\x S\to S$ is smooth which
implies that $\ps_{\be\al}:U_{\be\al}\to G$ is also smooth.
$(\ps_{\al\be})$ is a cocycle of transition functions and we
use it to glue a principal bundle with structure group $G$ over
$M$ which we call $(P,p,M,G)$. From its construction it is clear
that the associated bundle $P[S] = P\x_GS$ equals $(E,p,M,S)$.
\endremark

We have thus shown assertions \therosteritem1 and \therosteritem2 of 
the theorem. The two remaining assertions will be shown later in 
\nmb!{4.7}.
\qed\enddemo

\heading\totoc\nmb0{4}. The universal connection \endheading

\subheading{\nmb.{4.1}. The extension of the system} 
Let $(H,\et)$ be a complete strong system of vector fields on 
the bundle $E$, and let $(C=C(H),p_C,M)$ be the bundle of 
$H$-connections on $E$ described in \nmb!{2.6}. We consider the 
following fibered products
$$\CD
TC\x_{TM}H  @>pr_2>> H         @. \qquad @.   TC\x_MA  @>pr_2>>      A \\
@Vpr_1VV             @VV\uet V       @.      @Vpr_1VV               @VVp_AV\\
TC        @>>T(p_C)> TM,       @. \qquad @.   TC  @>>\pi_M\o T(p_C)> M.
\endCD$$
There is a smooth diffeomorphism, fibered over $TC$, between these 
two fibered products, whose description involves the properties of 
$C=C(H)$:
$$\gather
\ka: TC\x_{TM}H \to TC\x_MA\\
(X,h)\mapsto (X,h-\pi_C(X).\uet(h))\\
(X,a+\pi_C(X).Tp_C.X) \gets (X,a).
\endgather$$
Since $TC\x_MA=TC\x_Cp_C^*A\to C$ is a vector bundle, we may regard 
also $TC\x_{TM}H\to C$ as a vector bundle via $\ka$.

We consider now the fiber bundle $C(H)\x_ME=C\x_ME\to C$ with 
standard fiber $S$, and the extended system of vector fields on it 
which is given by the following diagram:
$$\CD
(TC\x_{TM}H)\x_C(C\x_ME)  @>\tilde\et>>   T(C\x_ME)\\
@|                                        @| \\
TC\x_{TM}(H\x_ME)      @>>TC\x_{TM}\et>   TC\x_{TM}TE.
\endCD$$
Since the vertical part of $\tilde \et$ is the same as that of $\et$ 
we see that the system  $(TC\x_{TM}H,\tilde\et)$ is again strong and 
complete if the system $H$ is it. Canonical atlases for the extended system 
are given by base extensions of the canonical atlases for $(H,\et)$.
The exact sequence in the sense of \nmb!{2.3} is here an exact 
sequence of vector bundles over $C$:
$$0 \to C\x_MA \to TC\x_MA @>pr_1>>  TC  \to 0.$$

\subheading{\nmb.{4.2}. The universal connection} We now consider the 
bundle of connections $C(TC\x_{TM}H) \to C$ for the extended system of 
vector fields, in the sense of \nmb!{2.6}. It is just the affine 
bundle of all splittings of the exact sequence of vector bundles over 
$C$ in the bottom line of the big diagram of \nmb!{4.1}. Thus we have
$$\align
C(TC\x_{TM}H) &= \{\si\in L(TC,TC\x_MA): pr_1\o \si =Id_{TC}\}\\
&\cong L(TC,C\x_MA),
\endalign$$
since this affine bundle has a canonical section, namely 
$(Id_{TC},0)$. This canonical section gives rise to a distinguished 
connection on the bundle $C\x_ME\to C$ which is called the 
\idx{\it universal connection} since it has the universal property 
described in lemma \nmb!{4.3} below. Its horizontal lift will be 
called 
$$C^{\text{univ}}: TC\x_C(C\x_ME)=TC\x_ME \to T(C\x_ME)=TC\x_{TM}TE.$$
By the general formula of \nmb!{2.6} we have (taking into account all 
isomorphisms):
$$\align
C^{\text{univ}}(X,e) &= \tilde\et(\ka\i(X,0),e)\\
&= (TC\x_{TM}\et) (X,(\pi_C(X).Tp_C.X),e)\\
&= (X,\et(\pi_C(X).Tp_C.X,e)).
\endalign$$
This coincides with the coordinate formula of \cite{Mo}.
The unversal connection itself is then given by 
$$\gather
\Ph^{\text{univ}}: T(C\x_ME) = TC\x_{TM}TE \to V(C\x_ME) = C\x_MTE\\
\Ph^{\text{univ}}(X,Y) = (\pi_C(X),Y-\et(\pi_C(X).Tp_C.X,\pi_E(Y))).
\endgather$$

\proclaim{\nmb.{4.3}. Lemma} The universal connection has the 
following universal property: Let $\si\in C^\infty(C(H))$ be a 
section describing a horizontal lift $C_\si$ of a $H$-connection on $E$.
Consider the extended section 
$\si\x_ME: E= M\x_ME\to C\x_ME$. Then the unversal connection 
$C^{\text{univ}}$ on $C\x_ME$ and the connection $C_\si$ are
$(\si\x_ME)$-related, i\. e\. the following diagram commutes:
$$\CD
TM\x_ME  @>C_\si>>  TE  \\
@VT\si\x_MEVV       @VVT(\si\x_ME)V\\
TC\x_ME  @>>C^{\text{univ}}>  T(C\x_ME).
\endCD$$
Likewise the vertical projection $\Ph_\si:TE\to VE$ and the vertical 
projection $\Ph^{\text{univ}}$ of the universal connection are 
$(\si\x_ME)$-related, i\. e\. the following diagram commutes:
$$\CD
TE=T(M\x_ME)            @>\Ph_\si>>             VE=V(M\x_ME)\\
@VT\si\x_{TM}TEVV                   @VV\si\x_M{\text{ins}}V\\
T(C\x_ME)=TC\x_{TM}TE   @>>\Ph^{\text{univ}}>   V(C\x_ME)=C\x_MTE.
\endCD$$
\endproclaim

\demo{Proof}
Check from the definitions that the diagrams commute.
\qed\enddemo

\subheading{\nmb.{4.4}. The universal holonomy group} 
Since the universal connection 
$C^{\text{univ}}$ respects the system $TC\x_{TM}H$ on $C\x_ME\to C$, 
and since this system is complete as noted in \nmb!{4.1},
$C^{\text{univ}}$ is a complete connection by theorem \nmb!{3.1}.

Now we choose $c_0\in C$ with $p_C(c_0)=x_0\in M$ and we identify 
again the standard fiber $S$ with $(C\x_ME)_{c_0}\cong E_{x_0}$. 
Then we can consider the holonomy group 
$\operatorname{Hol}^{\text{univ}}(c_0):=
     \operatorname{Hol}(C^{\text{univ}},c_0)$
within the group of all diffeomorphisms of the standard fiber $S$.
We may now apply the first half of the proof of theorem \nmb!{3.6} to 
the universal connection $C^{\text{univ}}$ on the bundle 
$C\x_ME\to C$. From step 3 of that proof it follows that the 
universal holonomy group $\operatorname{Hol}^{\text{univ}}(c_0)$ is a 
subgroup of the Lie group $G$ constructed there. The 
groups coincide, but we will not need this fact.

\proclaim{\nmb.{4.5}. Lemma} The parallel transport 
$\Pt^{\text{univ}}$ of the universal connection has the following 
universal property:

Let $\si\in C^\infty(C(H))$ be a 
section describing a horizontal lift $C_\si$ of a $H$-connection on $E$.
Let $c:[0,1]\to M$ be a (piecewise) smooth curve in $M$. Then the 
universal parallel transport $\Pt^{\text{univ}}$ 
and the parallel transport $\Pt^\si$ of the connection $C_\si$ are 
related by the following formulas:
$$\gather
\Pt^{\text{univ}}(\si\o c,t)\o(\si\x_ME) = (\si\x_ME)\Pt^\si(c,t)\\
pr_2\o\Pt^{\text{univ}}(\si\o c,t) = \Pt^\si(c,t)\o pr_2.
\endgather$$
\endproclaim

\demo{Proof}
We only have to show that for $u\in E_{c(0)}$ the following formula 
holds:
$$\Pt^{\text{univ}}(\si\o c,t,(\si(c(0)),u)) = (\si\x_ME)\Pt^\si(c,t,u)$$
Both curves cover the curve $\si\o c$ in $C$ and have the same 
initial value $(\si(c(0)),u)\in C\x_ME$. Moreover by lemma \nmb!{4.3} 
we have
$$\align
\Ph^{\text{univ}}\tfrac d{dt}(\si\x_ME)\Pt^\si(c,t,u)
&= \Ph^{\text{univ}}(T\si\x_{TM}TE)\tfrac d{dt} \Pt^\si(c,t,u)\\
&= (\si\x_MTE)\Ph_\si \tfrac d{dt} \Pt^\si(c,t,u) = 0. \qed
\endalign$$
\enddemo

\proclaim{\nmb.{4.6}. Lemma} Let $b:[0,1]\to C_{x}$ be a vertical 
(piecewise) smooth curve in $C$. Then the universal parallel 
transport along $b$ is just given by the affine structure of 
$C\to M$, i\. e\. we have 
$\Pt^{\text{univ}}(b,t,(b(0),e)) = ((b(t),e))$
for each $e\in E_x$
\endproclaim

\demo{Proof}
By the formula for $\Ph^{\text{univ}}$ in \nmb!{4.2} we have 
$$\align
\Ph^{\text{univ}} \tfrac d{dt} (b(t),0_e) 
     &= \Ph^{\text{univ}} (b'(t),0_e)\\
&= (b(t),0_e - \et(b(t).0_x,e)) = 0_{(b(t),e)}. \qed
\endalign$$
\enddemo

\subheading{\nmb.{4.7}. Rest of the proof of theorem \nmb!{3.6}}
We assume that we are again in the situation at the end of the proof.

\demo{Step 5} Lifting each $H$-connection to $P$. \newline
For this we have to compute the Christoffel symbols of $C_\ta$ for an 
arbitrary section $\ta\in C^\infty(C(H))$
with respect to the atlas of step 4. To do this directly is quite 
difficult since we have to differentiate the parallel transport with 
respect to the curve. Fortunately there is another way using the 
universal parallel transport. Let again $\Pt^\ta$ denote the parallel 
transport of $C_\ta$ and as above $\Pt=\Pt^\si$ that one of $C=C_\si$.
Let us identify 
$S\cong E_{x_0}\cong (C\x_mE)_{\si(x_0)}=\{\si(x_0)\}\x S$.
Let $c:[0,1]\to U_\al$ be a smooth curve. Then we have 
$$\align
\ps_\al&(\Pt^\ta(c,t)\ps_\al\i(c(0),s)) = \\ 
&=\left( c(t), \Pt(c_\al\i,1)
     \Pt((c_\al^{c(t)})\i,1)\Pt^\ta(c,t)\Pt(c_\al^{c(t)},1)
     \Pt(c_\al,1)s \right).
\endalign$$
Let now $b_0:[0,1]\to C_{c(0)}$ be a vertical smooth curve from 
$\si(c(0))$ to $\ta(c(0))$, and let $b_t:[0,1]\to C_{c(t)}$ be one 
from $\si(c(t))$ to $\ta(c(t))$. Using lemmas \nmb!{4.5} and \nmb!{4.6}
the last expression then gives
$$\align
\ps_\al&(\Pt^\ta(c,t)\ps_\al\i(c(0),s)) = \\ 
&=\Bigl( c(t), pr_2\Pt^{\text{univ}}(\si\o(c_\al^{c(t)}.c_\al)\i,2)
     \Pt^{\text{univ}}(b_t\i,1)\Pt^{\text{univ}}(\ta\o c,t)\\
&\qquad\quad  \Pt^{\text{univ}}(b_0,1)
     \Pt^{\text{univ}}(\si\o(c_\al^{c(t)}.c_\al),2)(\si(x_0),s)\Bigr)\\
&= (c(t),\ga(t).s),
\endalign$$
where $\ga(t)$ is a smooth curve in the holonomy group $G$ since we have
$\operatorname{Hol}^{\text{univ}}(\si(x_0))\subset G$ as remarked in 
\nmb!{4.4}. 
Now let $\Ga^\al_\ta\in\Om^1(U_\al,\X(S))$ 
be the Christoffel symbol of the connection $\Ph_\ta$ 
with respect to the chart $(U_\al,\ps_\al)$. 
 From the third proof of theorem \cite{Mi, 1.5} we have 
$$\ps_\al(\Pt^\ta(c,t)\ps_\al\i(c(0),s))=\left(c(t),\bar\ga(t,s)\right),$$
where $\bar\ga(t,s)$ is the integral curve through $s$ of the time 
dependent vector field $\Ga^\al_\ta(\frac d{dt}c(t))$ on $S$. But then we 
get 
$$\Ga^\al_\ta(\tfrac d{dt}c(t))(\bar\ga(t,s))= \tfrac d{dt}\bar\ga(t,s) = 
     \tfrac d{dt}(\ga(t).s) = (\tfrac d{dt}\ga(t)).s,$$
where $\tfrac d{dt}\ga(t)\o\ga(t)\i\in \g$. So $\Ga^\al_\ta$ takes values in 
the Lie sub algebra of fundamental vector fields for the action of 
$G$ on $S$. Theorem \cite{Mi, 2.5} shows that 
the connection $\Ph_\ta$ is induced from a principal connection 
$\om_\ta$ on $P$. 

Thus any $H$-connection on $E=P[S]$ is induced by a principal 
connection on $P$. By \nmb!{2.7} this also implies that the system 
$(H,\et)$ is induced from the system $TP/G$ of $G$-invariant 
projectable vector fields on $P$. 
\qed\enddemo


\Refs

\ref 
\by Almeida, R.; Molino, P.
\paper Suites d'Atiyah et feuilletages transversalement complets
\jour C. R. Acad. Sci. Paris
\vol 300, Ser. 1
\pages 13--15
\yr 1985
\endref


\ref 
\by Garcia P\'erez, P.
\paper Gauge algebras, curvature, and symplectic structure
\jour J. Diff. Geom.
\vol 12
\pages 209--227
\yr 1977
\endref

\ref 
\by Gauthier, Jean-Paul
\book Structure des syst\`emes non-lin\'eaires 
\publ \'Editions du CNRS
\publaddr Paris
\yr 1984
\endref

\ref 
\by Mackenzie, Kirill
\book Lie groupoids and Lie algebroids in differential geometry
\bookinfo London Math. Soc. Lecture Notes Ser. 124
\publ Cambridge Univ. Press
\publaddr Cambridge etc
\yr 1987
\endref

\ref   
\by Mangiarotti, L.; Modugno, M.   
\paper Connections and differential calculus on fibered manifolds. Applications to field theory.   
\paperinfo preprint 1989   
\endref

\ref   
\by Marathe, K. B.; Modugno, M.   
\paper Polynomial connections on affine bundles   
\paperinfo Pre\-print 1988   
\endref

\widestnumber\key{Mo}

\ref
\key Mi   
\by Michor, P. W.   
\paper Gauge theory for diffeomorphism groups   
\inbook Proceedings of the Conference on Differential Geometric Methods in Theoretical Physics, Como 1987, K. Bleuler and M. Werner (eds.)   
\publ Kluwer   
\publaddr Dordrecht   
\yr 1988   
\pages 345--371   
\endref 

\ref \key{}
\by Michor, P. W.
\book Gauge theory for fiber bundles 
\bookinfo Monographs and Textbooks in Physical Science 19 
\publ Bibliopolis
\publaddr Napoli
\yr 1991 
\endref

\ref   
\key Mo
\by Modugno, M.   
\paper An introduction to systems of connections   
\jour Sem Ist. Matem. Appl. Firenze   
\vol 7   
\yr 1986   
\pages 1--76   
\endref

\ref\key{}   
\by Modugno, M.   
\paper Systems of vector valued forms on a fibred manifold and applications to gauge theories,   
\inbook Lecture Notes in Math.   
\vol 1251   
\publ Springer-Verlag   
\yr 1987   
\endref 

\ref\key{}   
\by Modugno, M.   
\paper Linear overconnections   
\inbook Proceedings Journ. Relat. Toulouse   
\yr 1988   
\pages 155--170   
\endref

\ref\key{}   
\by Modugno, M.   
\paper Systems of connections and invariant Lagrangians   
\inbook Differential geometric methods in theoretical physics, Proc. XV. Conf. Clausthal 1986   
\publ World Scientific Publishing   
\publaddr Singapore   
\yr 1987   
\endref

\ref\key{}   
\by Modugno, M.   
\paper Jet involution and prolongation of connections   
\jour Casopis Pest. Math.   
\yr 1987   
\endref

\ref\key{}   
\by Modugno, M.; Ragionieri, R; Stefani, G.   
\paper Differential pseudoconnections and field theories   
\jour Ann. Inst. H. Poincar\'e   
\vol 34 (4)   
\pages 465--493   
\yr 1981   
\endref

\ref\key{} 
\by Palais, Richard S.
\paper A global formulation of the Lie theory of transformation groups 
\jour Mem. AMS
\vol 22
\yr 1957
\endref

\endRefs
\enddocument